\documentclass{article}

\usepackage{arxiv}

\usepackage[utf8]{inputenc} 
\usepackage[T1]{fontenc}    
\usepackage{hyperref}       
\usepackage{url}            
\usepackage{booktabs}       
\usepackage{amsfonts}       
\usepackage{nicefrac}       
\usepackage{microtype}      
\usepackage{lipsum}		
\usepackage{graphicx}

\usepackage{moreverb}

\usepackage{graphicx, subfigure}

\newcommand\BibTeX{{\rmfamily B\kern-.05em \textsc{i\kern-.025em b}\kern-.08em
T\kern-.1667em\lower.7ex\hbox{E}\kern-.125emX}}

\usepackage{graphicx}
\RequirePackage{fix-cm}

\def\beq{\begin{equation}}
\def\eeq{\end{equation}}
\def\baq{\begin{eqnarray}}
\def\eaq{\end{eqnarray}}
\def\bc{\begin{center}}
\def\ec{\end{center}}
\def\ds{\displaystyle}
\def\fai{\varphi}
\def\eps{\varepsilon}

\def\vtheta{\vartheta}
\def\gr{\gamma_{r}}
\def\gr1{\gamma_{r1}}

\def\ttt{\theta} 
\usepackage{hyperref}
 \usepackage{epsfig} 
 \usepackage{amssymb}
 \usepackage{amsmath}
\usepackage{multirow}
\usepackage{epsfig} 
\input{epsf}
\usepackage{graphicx}
\usepackage{graphicx, subfigure}
\usepackage[utf8]{inputenc}
\usepackage[english]{babel}

\def\beq{\begin{equation}}
\def\eeq{\end{equation}}
\def\baq{\begin{eqnarray}}
\def\eaq{\end{eqnarray}}
\def\bal{\begin{align} }
\def\eal{\end{align} }
\def\bc{\begin{center}}
\def\ec{\end{center}}
\def\ds{\displaystyle}
\def\fai{\varphi}
\def\eps{\varepsilon}

\def\vtheta{\vartheta}
\def\gr{\gamma_{r}}
\def\gr1{\gamma_{r1}}

\def\fai{\varphi}
\def\eps{\varepsilon}

\def\vtheta{\vartheta}
\def\gr{\gamma_{r}}
\def\gr1{\gamma_{r1}}

\def\fai{\varphi}
\def\eps{\varepsilon}

\def\vtheta{\vartheta}

\title{Systematic Review of  Newton-Schulz Iterations with Unified Factorizations : Integration in the Richardson Method  and  Application to Robust
Failure Detection in Electrical Networks}
\author{ Alexander Stotsky \\
   Department of Computer Science and Engineering \\
     Chalmers University of Technology  \\
     Gothenburg  SE - 412 96, Sweden  \\
	\texttt{alexander.stotsky@chalmers.se} \\
        \texttt{alexander.stotsky@telia.com}  }
\hypersetup{
pdftitle={A template for the arxiv style},
pdfsubject={q-bio.NC, q-bio.QM},
pdfauthor={David S.~Hippocampus, Elias D.~Striatum},
pdfkeywords={First keyword, Second keyword, More},
}
\begin{document}
\maketitle
\begin{abstract}
~Systematic overview of Newton-Schulz and Durand iterations with convergence analysis and
factorizations is presented
in the  chronological sequence
in unified framework. Practical recommendations for the choice of the order and factorizations
of the algorithms and
integration into Richardson iteration are given.  The simplest combination of Newton-Schulz and
Richardson iteration is applied to the parameter estimation problem associated with the failure detection
via evaluation of the frequency content of the signals in electrical network. The detection is performed
on real data for which the software failure was simulated, which resulted in the rank deficient information matrix.
Robust preconditioning for rank deficient matrices is proposed and the efficiency of the approach is demonstrated
by simulations via comparison with standard LU decomposition method.
\end{abstract}
\maketitle
\section{Introduction}
\noindent
Parameter estimation problems appear often in adaptive control, \cite{fom},\cite{io}
system identification,\cite{ljung}, signal processing, \cite{ljung}  and in
many other areas. Suppose that the process is described by
the following equation:
\beq
y_k = \fai_k^T \ttt_* \label{inyk}
\eeq
where $\fai_k$ is the regressor vector and $\ttt_*$
is the vector of unknown parameters to be estimated, $k=1,2,...$
\\ Introduction  of the model in the form
\beq
\hat{y}_k = \fai_k^T \ttt_k \label{inykm}
\eeq
and minimization of the following function associated with the mismatch between the model and the process
\beq
 E_k = \sum_{p = 1}^{p=k} (\hat{y}_{p}  - y_p )^2
 \label{inei}
 \eeq
yields the following restricted linear system of equations :
\begin{align}
    A_k \theta_k  &= b_k \label{inab} \\
    A_k &= \sum_{p = 1}^{p = k } \fai_{p}~ \fai_{p}^{T},~
    b_k = \sum_{p = 1}^{p = k} \fai_{p}~ y_p   \label{inps82}
\end{align}
which should be solved with respect to $\ttt_k$. Accurate and computationally
efficient solution of  (\ref{inab}) is the most significant part
of the estimation problem. The solution may benefit
from the properties of the matrix  $A_k$, which depend on the type of regressor
and the step number associated with the window size.
The matrix $A_k$ is SPD (Symmetric and Positive Definite)
matrix in many cases. For example, the matrix $A$ is SPD for the systems with harmonic regressor,
\cite{fom} -  \cite{st10}.
In general case  the equation (\ref{inab}) can be multiplied by $A^T_k$ (for any invertible matrix $A_k$) that transforms the system to the SPD case with
the Gram matrix, \cite{bjor}. Notice that the Gram matrix is often ill-conditioned matrix due to squaring
of the condition number \cite{bjor}.
\\ The matrix $A$ is SPD ill-conditioned (a) information matrix for systems with harmonic regressor
and short window sizes, \cite{cdc54},\cite{sto2019} (b) Gram matrix due to the squaring of the condition number in the least squares method, \cite{ljung}, \cite{bjor},  (c) mass matrix in finite element method, \cite{stic}  and mass matrix (lumped mass matrix) for  mechanical systems with singular perturbations, \cite{mostafa}, (d) state matrix for systems of linear equations, \cite{has}  and in many other applications.
Moreover, many processes in practice are singulary perturbed (have modes with different time-scales), \cite{kok}
and considered as stiff and ill-conditioned systems, which potentially extends application areas of  this model.
In other words the model (\ref{inab}) with SPD (and possibly ill-conditioned matrix $A_k$) covers almost all the cases in the parameter estimation applications.
\\ Ill-conditioning implies robustness problems (sensitivity to numerical calculations) and
imposes additional requirements on the accuracy of the solution of (\ref{inab})
(which are especially pronounced in finite-digit calculations) since small changes in $b_k$ due to measurement, truncation, accumulation, rounding and other errors result in significant changes in $\theta_k$.
In addition, ill-conditioning implies slow convergence of the iterative procedures.
\\ The accuracy requirements is the main motivation for application of the Richardson iteration which is driven by
the residual error and has filtering (averaging) property.
The residual error is smoothed and remains bounded for a sufficiently large step number
in this iteration, where the bound depends on inaccuracies, providing the best possible solution in finite digit calculations.
\\ Therefore convergence rate improvement and reduction of the computational complexity of
the Richardson iteration is one of the most important problems in the area. One of the ways to convergence rate
improvement is introduction of  the Newton-Schulz  and Durand algorithms in the Richardson iteration, \cite{sto2019}, \cite{dub} - \cite{st14}.
\\ This paper starts with short systematic overview (presented in unified tutorial fashion)
of  Newton-Schulz  and Durand iterations in Section~\ref{mia}.
Special attention is paid to practical introduction to power series factorizations
for reduction of the computational complexity in Section~\ref{fac}. The simplest integration of
Newton-Schulz  algorithms in the Richardson iteration is presented in Section~\ref{rich}
with application to practical problem in Section~\ref{gpm}.  The problem is associated with the failure
detection via estimation of the frequency content of the signals in electrical system
in the presence of missing data. The detection is performed on real signals
with simulated failure in the software. The efficiency of the approach is  demonstrated in  Section~\ref{gpm}.
The paper ends with brief conclusions in Section~\ref{con}.

\section{Overview of Iterative Matrix Inversion Algorithms: ~Mini Tutorial}
\label{mia}
\noindent
Brief overview of iterative matrix inversion algorithms
with convergence analysis is presented in this Section.
The overview is presented in the systematic way (from simple algorithms to
more sophisticated ones) in the  chronological sequence with convergence
analysis and starts with the second order Newton-Schulz  algorithms in Section~\ref{so}
followed by the description of high order Newton-Schulz algorithms in Section~\ref{hoi}.
Computationally efficient Durand iterations are described in Section~\ref{dui}, and
finally the combination for convergence rate improvement of high order Newton-Schulz and Durand algorithms
is presented in Section~\ref{conss}.
\subsection{Second Order Newton-Schulz Iterations}
\label{so}
\noindent
Consider the following iterations, \cite{shul} - \cite{beniz}:
\beq
G_k = G_{k-1} + F_{k-1} G_{k-1}, ~~ F_k = I - G_k A,     \label{dem1}
\eeq
where $I$ is the identity matrix, $k=0,1,2...$ .  Evaluation of the error in the first step $F_1$ yields :
\begin{align}
F_1 &= I - (G_0  + F_0~G_0)~A = I - (I + F_0) ~G_0 ~A    \nonumber \\
     &= I - (I + F_0)~(I-F_0) = I - I^2 + F_0^2 = F_0^2   \nonumber  \\
  F_2 &= F_1^2 = F_0^4  \nonumber
\end{align}
\noindent
which gives the error model:
\beq
  F_k = F_0^{2^k}   \label{demerm1}
\eeq
Algorithm (\ref{dem1}) which is also known as  Newton-Schulz-Hotelling matrix iteration, \cite{hot}
estimates the inverse of the matrix $A$ via minimization of the estimation
error $F_k$ provided that the spectral radius
of initial matrix
\beq
  \rho(F_0) = \rho(I - G_0 A) < 1  \label{demcon1}
\eeq
is less than one, where $G_0$ is preconditioner (initial guess of $A^{-1}$). The proof of convergence is presented in  \cite{shul} - \cite{beniz}, where
the algorithm (\ref{dem1}) is written in the following form:
$G_k = G_{k-1} + G_{k-1} F_{k-1}$, where $F_k = I - A G_k $.
\\ The order of the iterative algorithm can be increased for convergence rate improvement
by introduction of the power series in the algorithm (\ref{dem1}).
\subsection{High Order Newton-Schulz Iterations}
\label{hoi}
\noindent
High order algorithms (which are also known as hyperpower iterative algorithms, which are widely used
for calculation of generalized inverses)
can be presented in the following form for $n=2,3,...$, \cite{isa} - \cite{stekel} :
\begin{align}
G_k &= \sum_{j=0}^{n-1} F^j_{k-1} G_{k-1} = (I + F_{k-1} + F^2_{k-1} + ... + F^{n-1}_{k-1})~ G_{k-1}
\label{ik1}
\end{align}
or in the form presented in \cite{st14}, \cite{petru} :
\begin{align}
G_k &= \sum_{j=0}^{n-1} F^j_{k-1} G_{k-1} =  G_{k-1} + F_{k-1}~\sum_{j=0}^{n-2} F^j_{k-1} G_{k-1}
\label{pet1}
\end{align}
The error $F_k$ in (\ref{ik1}) can be written in the following form:
\begin{align}
F_k &= I - G_k A = I -(I + F_{k-1} + F^2_{k-1} + ... + F^{n-1}_{k-1}) G_{k-1} A \nonumber \\
 &= I -(I + F_{k-1} + F^2_{k-1} + ... + F^{n-1}_{k-1}) (I - F_{k-1}) = F_{k-1}^n \nonumber \\
  F_k &= F_0^{n^k}   \label{isa2}
\end{align}
Notice that the algorithm (\ref{ik1}) is written in following form in \cite{isa} - \cite{stekel} (and it is  widely used in this form):
\begin{align}
G_k &= G_{k-1} \sum_{j=0}^{n-1} F^j_{k-1} = G_{k-1} (I + F_{k-1} + F^2_{k-1} + ... + F^{n-1}_{k-1})
\label{ik2}
\end{align}
with $F_k = I - A G_k $.
Notice that the error $F_k = I - G_k A$ is better aligned with Richardson
iteration, see Section~\ref{rich} and therefore is considered further in this paper.
\subsection{Durand Iterations}
\label{dui}
\noindent
Computationally efficient matrix inversion method
described by Durand, \cite{durand}, \cite{bre}
can be derived from the relation associated with the initial error in Newton-Schulz approach,~ $A^{-1} = F_0 A^{-1} + G_0$ by substitution of $G_k$ instead of $A^{-1}$ as follows :

\begin{align}
G_k &= F_0 ~G_{k-1} + G_0 \label{duran1} \\
G_k - A^{-1} &= F_0~ (G_{k-1} - A^{-1}) \label{ermoddr} \\
F_k &= I - G_k A = F^{k+1}_0  \label{ermoddr1}
\end{align}
The error models (\ref{ermoddr}), (\ref{ermoddr1}) are valid for the algorithm
(\ref{duran1}).
Higher order Durand iterations are
discussed in \cite{sto2019}, \cite{bre}.
\\ Notice that the convergence rate of Newton-Schulz Iterations is significantly higher
than the  convergence rate of Durand iterations. However, Durand iteration (\ref{duran1})
requires only one matrix product in each step.
\\ Notice also that the convergence performance can be enhanced by introduction
of the initial power series as preconditioning $G_0$
(which is calculated only once) in Newton-Schulz  and Durand iterations,
\cite{sto2019}, \cite{or}.
The order of initial power series can be chosen sufficiently high, which essentially
improves the convergence rate.  Computational complexity of initial power
series can be reduced via factorizations,
see Section~\ref{fac} for details.
\\ Convergence rate can also be improved via combination of two
Newton-Schulz iterations, which also cover Durand iterations as special case, see
next Section~\ref{conss}.

\subsection{Convergence Rate Improvement via Combination of High Order Newton-Schulz Iterations}
\label{conss}
\noindent
Consider the following combination of two Newton-Schulz iterations associated with
two-step iterative method, \cite{ifac2020} :
\begin{align}
\Gamma_k &= I - L_{k-1} A \label{gamk11} \\
L_k &=  \{ \sum_{j=0}^{n-1} \Gamma^j_k \} ~ L_{k-1} \label{lkk1} \\
\Gamma^n_k &=  I - L_k A \label{gamk21}  \\
G_k &= \underbrace{L_k}_{\begin{subarray}{l}\text{Newton-Schulz }\\
    \text{Iteration }\end{subarray}} + \underbrace{\Gamma^n_k}_{\begin{subarray}{l}\text{High Order }\\
    \text{Convergence} \\
    \text{Accelerator} \end{subarray}}   ~  ~ \underbrace{ \{ \sum_{j=0}^{n-1}  F^j_{k-1} \} ~ G_{k-1}}_
 {\begin{subarray}{l}\text{Newton-Schulz }\\
    \text{Iteration }\end{subarray}}
 \label{gkgen21} \\
 F_k &= I -  G_k A = \Gamma^n_k ~ F^n_{k-1}      \label{fks1} \\
F_k &= F_0^{\ds ~  k ~ n^{k+1} + ~ n^k },~\text{for}~n > 1  \label{ermoex1}
\end{align}
where $ \ds L_0 = \sum_{j=0}^{n-1} \Gamma^j_0 ~ S^{-1} $,~  $ \Gamma_0 = F_0 = I - S^{-1} A$, where
$\ds G_0 = S^{-1}$ satisfies (\ref{demcon1}).
Notice that Durand iteration (\ref{duran1}) can be derived from (\ref{gamk11}) - (\ref{gkgen21}) with $n=1$.
The error model (\ref{fks1}) can be  proved via multiplication of (\ref{gkgen21}) by $A$ and
explicit evaluation of the error $F_k$ (similar to Section~\ref{hoi}, see also \cite{isa}, \cite{ifac2020}) as follows:
\begin{align}
F_k &= I - G_k A =  \underbrace{I - L_k A}_{\Gamma^n_k} - \Gamma^n_k  \sum_{j=0}^{n-1}  F^j_{k-1} \underbrace{G_{k-1} A}_{I - F_{k-1}} = \Gamma^n_k ~ F^n_{k-1} \label{prof1}
\end{align}
The error model (\ref{ermoex1}) is obtained via explicit evaluation of (\ref{prof1}).
 The algorithm (\ref{gamk11}) - (\ref{ermoex1}) has two independent and equally complex computational parts: the first part is associated with calculations of $\Gamma_k$, $L_k$ and $\Gamma^n_k$, where $\Gamma_k = \Gamma^n_{k-1}$, whereas the term  $\ds  \sum_{j=0}^{n-1}  F^j_{k-1}  G_{k-1}$ is calculated in the second part. This algorithm is ideally suited for parallel implementation.
\\ \\ Extended version of this survey is presented in
Stotsky A., Recursive Versus Nonrecursive Richardson Algorithms:
 Systematic Overview, Unified Frameworks and  Application to Electric Grid Power Quality Monitoring,
 Automatika, vol. 63, N2, pp. 328-337, 2022, https://doi.org/10.1080/00051144.2022.2039989,
\cite{sto2021}.
\\ \\
Further developments can be found in Stotsky A.,Simultaneous Frequency and Amplitude Estimation
for Grid Quality Monitoring : New Partitioning with Memory Based Newton-Schulz Corrections,
IFAC PapersOnLine 55-9, pp. 42-47, 2022, https://www.sciencedirect.com/science/article/pii/S2405896322003937  ,\cite{ifac2022}.
\\

\section{Reduction of Computational Complexity via Unified Factorizations}
\label{fac}
\noindent
The following high order Newton-Schulz power series can be factorized
for composite orders $n$ as follows, \cite{ifac2020} :
\begin{align}
G_k &= \{ \sum_{j=0}^{n-1} F_{k-1}^j \} ~ G_{k-1}  = \{ \sum_{j=0}^{w-1} F_{k-1}^{(p+1) j} \} ~
\{ \sum_{d=0}^{p} F_{k-1}^d \} ~ G_{k-1}    \label{ztab1} \\
   F_{k-1}^{p+1} &= I -  \{ \sum_{d=0}^{p} F_{k-1}^d \} ~ G_{k-1}~A  \nonumber \\
   &=  I -  \{ \sum_{d=0}^{p} F_{k-1}^d \} ~ \{ I - F_{k-1} \}, ~ F_{k-1} = I - G_{k-1} A       \label{ua1} \\
    n &= w ~ ( p + 1), ~~ p = 0,1,2, ..., ~ w = 1,2,3, ...   \label{hape}
\end{align}
Taking into account that every prime number is followed after the composite
number the factorization for the prime orders $n + 1$  can be performed as follows:
\begin{align}
  G_k  &= (I + F_{k-1} \{ \sum_{j=0}^{w-1} F_{k-1}^{(p+1) j} \} ~ \{ \sum_{d=0}^{p} F_{k-1}^d \}) ~ G_{k-1}  \label{ztab2}  \\
   &= \{ \sum_{j=0}^{(n+1) - 1} F_{k-1}^j \} ~ G_{k-1}  \nonumber
\end{align}
Similar factorizations can be found in \cite{solei1}. The following efficiency index
introduced in \cite{ehr} is widely used for quantification of the performance of
different factorizations:
\beq
 EI = \ds  n^{1 / n_p} \label{ei}
\eeq
where $n$ is the order of the algorithm and $n_p$ is the number of matrix products per each cycle.
For illustration of the unified factorization approach few examples with increasing EI are presented below.
\subsection{Case of $n = 8 = w ~ ( p + 1)$,~ $p = 3$,~$w = 2$  }
\begin{align}
G_k &= \{ \sum_{j=0}^{1} F_{k-1}^{ 4 j} \} ~ \{ \sum_{d=0}^{3} F_{k-1}^d \} ~ G_{k-1} \label{o8}  \\
 &=   (I +  F_{k-1}^4) (I +  F_{k-1}^2) (I +  F_{k-1}) ~ G_{k-1} \nonumber
\end{align}
which requires $6$ matrix products with $EI =  1.4142$ .

\subsection{Case of $n = 9 = w ~ ( p + 1)$,~ $p = 2$,~$w = 3$ }
\begin{align}
G_k &= \{ \sum_{j=0}^{2} F_{k-1}^{ 3 j} \} ~ \{ \sum_{d=0}^{2} F_{k-1}^d \} ~ G_{k-1}  \nonumber \\
&=   (I +  F_{k-1}^3 +  F_{k-1}^6) (I +  F_{k-1}  +  F_{k-1}^2 ) ~ G_{k-1} \nonumber \\
&=    (I + (I + F_{k-1}^4) (I + F_{k-1}^2) (F_{k-1} + F_{k-1}^2) ) ~ G_{k-1} \nonumber
\end{align}
which requires $6$ matrix products with  $EI =   1.4422$ .
\subsection{Case of $n = 10 = w ~ ( p + 1)$,~ $p = 1$,~$w = 5$ }
\begin{align}
G_k &= \{ \sum_{j=0}^{4} F_{k-1}^{ 2 j} \} ~ \{ \sum_{d=0}^{1} F_{k-1}^d \} ~ G_{k-1}  \nonumber \\
&=   (I  + (F_{k-1}^2 + F_{k-1}^4) (I + F_{k-1}^4) ) (I + F_{k-1})  ) ~ G_{k-1} \nonumber
\end{align}
which requires $6$ matrix products with $EI =1.4678$ .

\subsection{Case of $n = 11 = w ~ ( p + 1) + 1$,~ $p = 1$,~$w = 5$ }
\begin{align}
G_k &= (I +  F_{k-1} \{ \sum_{j=0}^{4} F_{k-1}^{ 2 j} \} ~ \{ \sum_{d=0}^{1} F_{k-1}^d \}) ~ G_{k-1}
\label{o11} \\
  &=   (I + F_{k-1} (I  + (F_{k-1}^2 + F_{k-1}^4) (I + F_{k-1}^4) ) (I + F_{k-1})  )  ~ G_{k-1} \nonumber \\
   &=   (I + (I  + (F_{k-1}^2 + F_{k-1}^4) (I + F_{k-1}^4) ) ( F_{k-1} + F^2_{k-1})  )  ~ G_{k-1} \nonumber
\end{align}
which requires  $6$ matrix products,  $EI = 1.4913$,  \cite{stan11}.
\\ Notice that the idea of factorization is associated with Newton-Schulz method.
Therefore Newton-Schulz  algorithms of low orders (order two or three)
being iterated for a number of steps can be applied instead of factorized Newton-Schulz  iterations of higher orders for the sake of robustness and efficiency.
\\
Notice that three steps of the
second order Newton-Schulz iteration are equivalent to one step
of eighth order iteration with $6$ matrix products, see (\ref{o8}).
However, application of the eleventh order algorithm, see (\ref{o11}) and
\cite{stan11} requires also $6$ mmm and provides faster convergence
than three steps of the second order Newton-Schulz iteration.
This algorithm has also the highest efficiency index, $EI = 1.4913$
among the algorithms presented above.
Therefore the proper choice of the order and factorization that is made for each particular
application should represent the trade-off between the robustness and convergence rate.
\section{Richardson Iterations}
\label{rich}
\noindent
Many practical  problems can be associated with
restricted linear system of equations:
\beq
   A \theta_* = b \label{riab}
\eeq
which should be solved for $\theta_*$. Richardson iteration is the most promising solution
to this problem which provides the best possible accuracy in
finite digit calculations.
Despite the fact that the inverse of matrix $A$ is not required in order to find $\theta_*$
the approximate inverse essentially improves convergence rate of the
Richardson iteration, \cite{sto2019}, \cite{dub} - \cite{st14}.  The simplest combination
of matrix inversion techniques and Richardson algorithm can be presented as follows :
\begin{align}
 \theta_k &=  \theta_{k-1} - G_k ~(A \theta_{k-1}  - b) = F_{k-1}  \theta_{k-1} + G_k b  \label{riri}
\end{align}
where  $\theta_k$  is the estimate of $\theta_*$, ~$G_k$ is associated with the estimate of the inverse of $A$ and $F_{k}$ is the inversion error,
see  Section~\ref{mia}.
Iteration  (\ref{riri}) can be implemented using matrix vector products only.
\\ Notice that definition of the inversion error in the form $F_k = I - G_k A$ is more convenient than
the following definition  $F_k = I - A G_k $ since the error $F_k = I - G_k A$  is better aligned with
Richardson iteration providing the same convergence conditions.
The convergence analysis of Richardson iteration with the error  $F_k = I - A G_k $
can be found for example in \cite{sri1}.
\\ Any matrix inversion algorithms presented in Section~\ref{mia} with factorizations described
in Section~\ref{fac} can be applied in (\ref{riri}).
 Error models for different choices of $G_k$ are presented in \cite{sto2019}.
Since $G_k$ is getting closer to $A^{-1}$ in each step the convergence rate of (\ref{riri}) is essentially improved.
Updating of $G_k$  can even be stopped after several steps for reduction of the computational burden.
 \begin{figure*}[ht]
     \begin{center}
        \subfigure[  ] {%
        \label{suba}
\includegraphics[height=60mm]{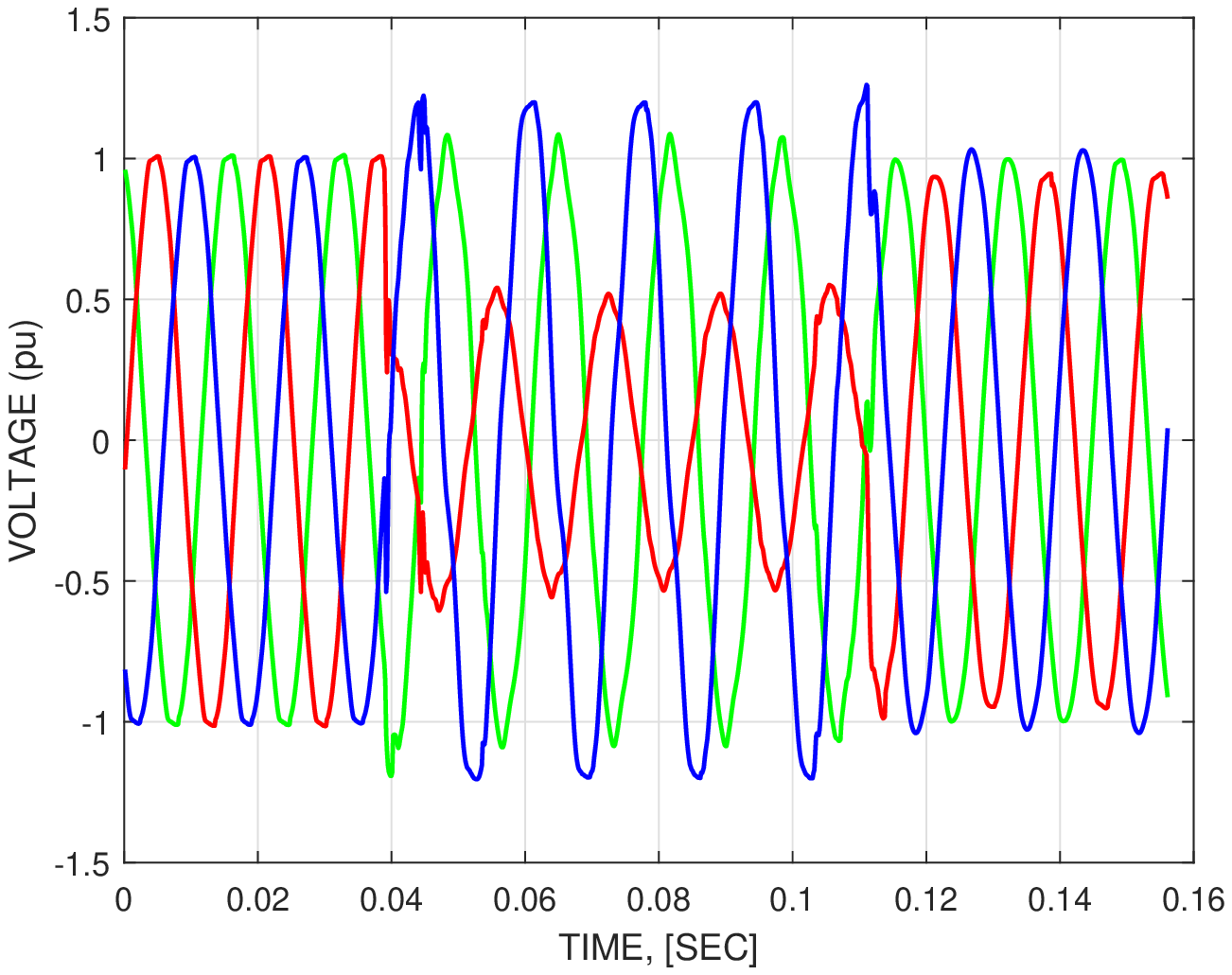} }%
 \subfigure[  ]{%
\label{subc}
\includegraphics[height=60mm]{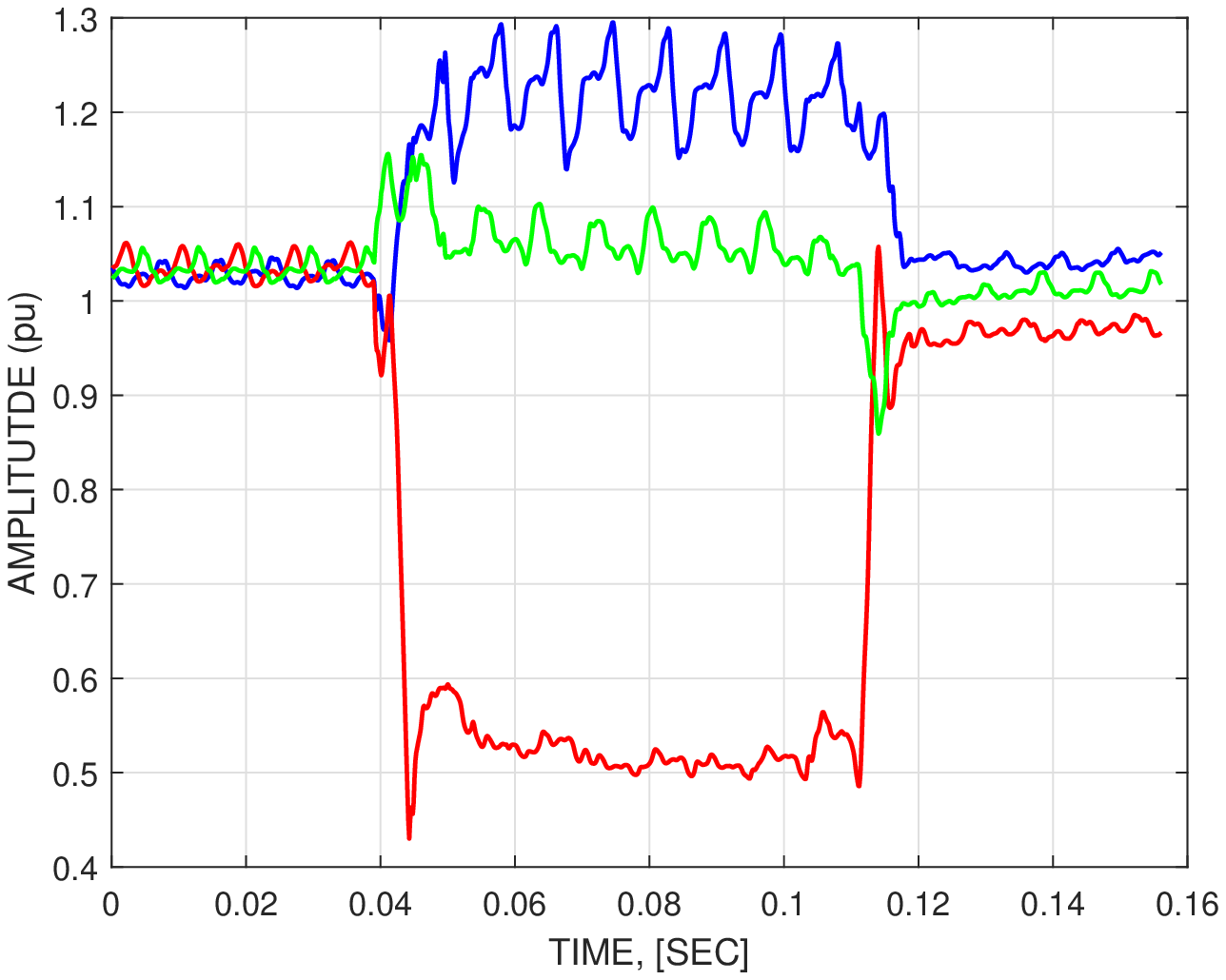}}%
\end{center}
\caption{
The voltages which correspond to the sag and swell events
in phases 'a','b' and 'c' are plotted with green, red and blue lines respectively
in Figure~\ref{suba}.
The voltages are associated with the  failure event with ID $2911$.
The data is obtained from the DOE/EPRI National Database of Power System Events, \cite{doe}.
The cause of this event is unknown.
Amplitudes of the first harmonic (with the corresponding colours)
which show delectability of the sag and swell events
for the signals in Figure~\ref{suba} are plotted in  Figure~\ref{subc}.
The detection is performed with rank deficient information matrix.}
\label{Fig2}
\end{figure*}

 \begin{figure}
\centerline{\psfig{figure=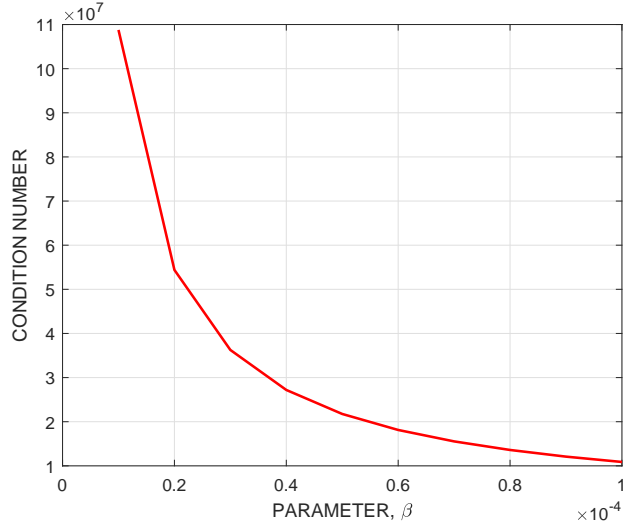,height=72mm}}
\begin{center}
\caption{\small{
The condition number of the matrix  $A_{rk}$ in (\ref{ar}) as a function of the regularization parameter $\beta >0$ is plotted.
For a small $\beta$ the matrix $A_{rk}$
is ill-conditioned matrix for which the spectral radius
$\ds  \rho(I - S_k^{-1} A_{rk})$, where $\ds S_k^{-1} = I / \alpha_k$ with $\alpha_k = \| A_{rk}\|_\infty / 2 + \eps$ is very close to one. }}
\label{fig1}
\end{center}
\end{figure}

 \section{Estimation of the Frequency Content of Electrical Signals for Robust Failure Detection}
\label{gpm}
\noindent
Suppose that the measured electrical  signal (voltage or current) $y_k$  can be presented
in the following form
\beq
  y_k = \fai_k^T \vtheta + \xi_k \label{ps1}
\eeq
where $\vtheta$ is the vector of
unknown constant parameters and   $\fai_k$ is the  harmonic regressor
presented in the following form:
\baq
      & & \fai_k^{T} =  [ \cos( q_0 k ) ~ \sin(q_0 k )  ~ \cos( 2 q_0 k )
     \nonumber \\
      & & ~ \sin( 2 q_0 k )~ ...~ \cos( m q_0 k ) ~ \sin( m q_0 k ) ]  \label{f1}
\eaq
where $q_0$ is the fundamental frequency of network (for example $q_0 = 50$ Hertz or $q_0 = 60$ Hertz),
$m$ is the number of harmonics, and $\xi_k$ is a zero mean
white Gaussian noise, $k = 1,2,...$ is the step number.
\\ The  model of the signal (\ref{ps1}) with adjustable parameters $ \theta_{*k}  $
is presented in the following form:
\beq
  \hat{y}_k =  \fai_k^{T}  \theta_{*k}  \label{ps4}
\eeq
The signal $y_k$ is approximated by the model $\hat{y}_k$  in the least squares sense in each step $k$ of moving window of a size $s$.
\\
The estimation algorithm is based on minimization of the following error $E_k$;
\beq
 E_k = \sum_{p = k -(s-1)}^{p=k} (\hat{y}_{p}  - y_p )^2
 \label{ei}
 \eeq
for a fixed step $k$, where $k \ge s$.
\\
The least squares solution for estimation of the parameter vector $\theta_{*k}$, which should be calculated in each step can be written as follows:
\baq
    A_k \theta_{*k} &=& b_k \label{ps8} \\
   A_k &=& \sum_{p = k -(s-1)}^{p=k} \fai_{p}~ \fai_{p}^{T}   \label{ps81} \\
    b_k &=& \sum_{p = k -(s-1)}^{p=k} \fai_{p}~ y_p   \label{ps82}
\eaq
\subsection{Missing Input Data \& Rank Deficiency of Information Matrix}
\noindent
The matrix $ A_k $ is symmetric and positive definite information matrix for
systems with harmonic regressor for a sufficiently large window size $s$.
However, the software, functional or system failures may occur, which usually result
in declaration of zero values in the information matrix  $A_k$. Such failures may result
in rank deficiency of the  information matrix $A_k$ in each step.
Nevertheless, the methods described above
can be still applied in this case.

\subsection{Fault Detection on Real Data }
\noindent
Power quality monitoring and detection of sag and swell events is
considered in the presence of system failures.
Monitoring is performed via estimation of the frequency content of the voltage signals in three phase
system. Real data obtained from the DOE/EPRI National Database of Power System Events, \cite{doe}
is used for power quality monitoring.
The event (Event ID $2911$) which resulted in voltage sag and swell was considered.
The cause of this event is unknown and the event is considered as temporary fault.
Measured voltage signals which correspond to the sag and swell
in phases 'a','b' and 'c' are plotted with green, red and blue lines respectively
in Figure~\ref{suba}.
\\
It is assumed that the voltage signals can be described by the model (\ref{ps1})
with harmonic regressor (\ref{f1}) and fundamental frequency of $60$ Hertz with four higher harmonics and unknown parameter vector $\vtheta$. The model of the signals is given by  (\ref{ps4}) with the vector of the parameters $\theta_{*k}$
which estimates unknown vector $\vtheta$. The estimation problem is reduced to the restricted linear system of equations (\ref{ps8}) which should be solved w.r.t. $\theta_{*k}$  in each step.
The $10~\times~10$ matrix $A_k$ is the SPD  matrix
in each step (which is verified by simulations).
\\ Suppose that a functional or system failure occurred, which resulted in three zero columns
(columns 3,4 and 5) of the information matrix  $A_k$ in each step.
\\ In addition, suppose that three components of the vector $b_k$ (the components 3,4 and 5)
are zeros in each step due to similar failure. Such failures result in rank deficient
information matrix $A_k$ in each step. The rank of the information matrix is equal to seven
in this case.  Nevertheless, the methods described above can be applied in this case.

\subsection{Pre-Conditioning for Rank Deficient Information Matrix}
\label{splpre}
\noindent
In order to apply the methods described above the matrix $S_k^{-1}$ should be chosen
such that the  spectral radius of the matrix $\ds I - S_k^{-1} A_k$   is less than one in
each step, $\ds  \rho(I - S_k^{-1} A_k) < 1 $. For SPD information matrices
the pre-conditioner
can be  chosen as  $\ds S_k^{-1} = I / \alpha_k$ with $\alpha_k = \| A_k\|_\infty / 2 + \eps$,
where $ \| \cdot \|_\infty $ is the maximum row sum matrix norm, $\eps > 0$,  \cite{cdc54}, \cite{ifac2017}.
In the case of software or system failures the information matrix may become rank deficient
and special preconditioning should be applied.
Consider the following transformation of the system (\ref{ps8}) and regularization of the matrix $A_k$:
\begin{align}
    A_{rk} \theta_{*k} &= A^T_k b_k \label{ar0} \\
     A_{rk} &= \beta I + A^T_k A_k  \label{ar}
\end{align}
where $\beta$ is a small positive number, and the matrix $A_{rk}$ is SPD matrix
for which the preconditioning described above can be applied.
The choice of the parameter $\beta$ represents the trade-off between high condition number (for very small
$\beta$) and the performance of regularization and hence the performance of the event
detection (for large $\beta$).
The condition number of  $A_{rk}$  as a function of the regularization parameter $\beta$ is shown in
Figure~\ref{fig1}. For a small $\beta$ which allows high performance detection the matrix $A_{rk}$
becomes ill-conditioned matrix  (which is sensitive to numerical operations) for which the spectral radius
$\ds  \rho(I - S_k^{-1} A_{rk})$ is very close to one.
Application of the stepwise splitting method,
where the matrix $A$ is split recursively in each step provides robust solution for ill-conditioned case, \cite{ifac2017}.
\\ Notice that the model (\ref{ar0}), (\ref{ar})
provides robust solution for parameter estimation problems
in the case of insufficient excitation
and singular information matrices and can be applied instead of classical self-tuning  adaptive and estimation
algorithms in adaptive control and system identification, \cite{fom} - \cite{ljung}.

\subsection{Simulation Results}
\label{simu}
\noindent
Amplitudes of the first harmonic of the measured signals
recovered by the Richardson parameter estimation algorithm  (\ref{riri}) with
the second order matrix inversion method (\ref{dem1})
are presented in Figure~\ref{subc}.
 The Figure shows that the sag and swell events are detectable using the approach described above
even in the presence of system failures, which result in rank deficiency of information matrix.
The estimation accuracy was evaluated  as the average sum of squares
of the parameter errors for Richardson algorithm with respect to
parameter vector calculated without missing data. The parameter vector calculated
with Richardson algorithm is quite close to the parameter vector calculated without missing data.
\\ In addition, LU decomposition method was also applied for parameter calculation in (\ref{ps8})
with regularized matrix $A_{rk}$, which is close to singular matrix (extremely ill-conditioned)
for sufficiently small $\beta$.  LU decomposition method is very sensitive to numerical errors
for extremely ill-conditioned matrices, whereas Richardson algorithm  delivers robust solution.
Therefore Richardson algorithms are recommended for robust detection of the failure events in the presence of missing data.

\section{CONCLUSIONS}
\label{con}
\noindent
The practical guide for solution of the parameter estimation problem
described by  restricted linear system of equations  was presented in the Richardson framework,
which provides the best possible solution in finite digit
calculations.
Newton-Schulz  and Durand iterations with unified factorizations, which enhance the performance of the
Richardson method   are presented in
the form of tutorial (in historical framework) with the convergence analysis.
\\ The combination of Richardson method with second order Newton-Schulz  approach
was tested in the failure detection problem on real signals from the  electrical network.
In addition, the software or system failure was simulated which resulted in rank deficiency of
the information matrix. New preconditioning method was proposed for rank deficient matrices
in the problem of estimation of the frequency contents of electrical signals.
It was shown that the Richardson method is applicable
in the failure detection problems with rank deficient information matrices.



\end{document}